\documentclass[11pt]{article}

\usepackage{amssymb}
\usepackage{amsmath}
\usepackage{stmaryrd}

\usepackage[latin1]{inputenc}

\newcommand {\limp}        {\mathbin{\Rightarrow}}
\newcommand {\leqv}        {\mathbin{\Leftrightarrow}}
\newcommand {\lconj}       {\mathbin{\wedge}}
\newcommand {\ldisj}       {\mathbin{\vee}}
\newcommand {\lneg}        {\mathop{\neg}}

\newcommand{\sat}     {\Vdash}
\newcommand{\ent}     {\vDash}
\newcommand{\der}     {\vdash}

\newcommand{\termvalue}[1] {\lbrack\!\lbrack #1 \rbrack\!\rbrack}

\newcommand{\bnfd}                               {\talloblong}

\newcommand{\cS}                                 {\mathcal{S}}
\newcommand{\cR}                                 {\mathcal{R}}
\newcommand{\ucS}                                {\cup\mathcal{S}}

\newcommand{\tuple}[1]                           {{(#1)}}
\newcommand{\ext}[1]                             {{|#1|}}
\newcommand{\prob}[1]                            {{(\int\!#1)}}

\newcommand{\qif}[3]                             {{(#1\rhd\,#2;\,#3)}}
\newcommand{\qiposs}[2]                          {{([#1]\lposs\;#2)}}
\newcommand{\bra}[1]                             {{\langle#1|}}
\newcommand{\ket}[1]                             {{|#1\rangle}}
\newcommand{\bkt}[2]                             {{\langle#1|#2\rangle}}

\newcommand{\qmols}                              {{\textrm{qmols}}}

\newcommand{\nats}                               {{\mathbb N}}

\newcommand{\reals}                              {{\mathbb R}}
\newcommand{\comps}                              {{\mathbb C}}

\newcommand{\qimp}                              {\mathbin{\sqsupset}}
\newcommand{\qconj}                             {\mathbin{\sqcap}}
\newcommand{\qdisj}                             {\mathbin{\sqcup}}
\newcommand{\qneg}                              {\mathop{\boxminus}}
\newcommand{\qeqv}                              {\mathbin{\equiv}}

\newcommand{\qAtom}                              {{\mathsf{qAtom}}}
\newcommand{\qB}                                 {{\mathsf{qB}}}
\newcommand{\qb}                                 {{\mathsf{qb}}}
\newcommand{\verum}                              {{\top}}
\newcommand{\falsum}                             {{\bot}}
\newcommand{\lnec}                               {{\Box}}
\newcommand{\lposs}                              {{\Diamond}}

\newcommand{\re}                                 {{\mathrm{Re}}}
\newcommand{\im}                                 {{\mathrm{Im}}}

\newcommand{\bfw}                                {{\mathbf w}}
\newcommand{\cB}                                 {{\mathcal{B}}}

\newcommand{\Prob}                               {{\mathcal{P}}}
\newcommand{\Hil}                                {{\mathcal{H}}}
\newcommand{\conjg}[1]                           {{\overline{#1}}}
\newcommand{\norm}[1]                            {{||#1||}}

\newcommand{\supp}                               {{\mathrm{supp}}}

\newcommand{\cata}                               {{\textbf{cat-alive}}}
\newcommand{\cati}                               {{\textbf{cat-in-box}}}
\newcommand{\catm}                               {{\textbf{cat-moving}}}

\begin{document}

\title{Weakly complete axiomatization of exogenous
quantum propositional logic}

\author{P.~Mateus and A.~Sernadas\\
CLC, Department of Mathematics, IST\\
\{pmat,acs\}@math.ist.utl.pt}

\maketitle

\begin{abstract}
A weakly complete finitary axiomatization for EQPL (exogenous quantum propositional logic) is presented. The proof is carried out using a non trivial extension of the Fagin-Halpern-Megiddo technique together with three Henkin style completions.
\end{abstract}

\section{Introduction}

A new logic (EQPL -- exogenous quantum propositional logic) was proposed in \cite{pmat:acs:04a,pmat:acs:04b} for modeling and reasoning about quantum systems, embodying all that is stated in the postulates of quantum physics (as presented, for instance, in \cite{coh:diu:lal:77,nie:chu:00}). The logic was designed
from the semantics upwards, starting with the key idea of adopting
superpositions of classical models as the models of the proposed quantum logic. 

This novel approach to quantum reasoning is quite different
from the traditional approach \cite{fou:99,chi:giu:gre:04} to the
problem that, as initially proposed by Birkhoff and von Neumann
\cite{bir:neu:36}, focuses on the lattice of closed subspaces of a
Hilbert space. Our exogenous semantics approach has the advantage of closely guiding the
design of the language around the underlying concepts of quantum physics
while keeping the classical connectives and was inspired by the possible
worlds approach originally proposed by Kripke \cite{kri:63} for modal
logic. It is also akin to the society semantics introduced in
\cite{car:mar:99} for many-valued logic and to the possible translations
semantics proposed in \cite{car:00} for paraconsistent logic. The
possible worlds approach was also used in
\cite{nil:86,nil:93,bac:90a,bac:90b,fag:hal:meg:90,aba:hal:94} for probabilistic
logic. Our semantics of quantum logic, although inspired by modal logic,
is also completely different from the alternative Kripke semantics given
to traditional quantum logics (as first proposed in \cite{dis:72}) still
closely related to the lattice-oriented operations. For other examples of logics based on the exogenous semantics approach see \cite{pmat:acs:css:05}.

Contrarily to traditional quantum logics that replace the classical
connectives by new connectives inspired by the lattice-oriented
operations, by adopting superpositions of classical models as the models
of the quantum logic, we are led to a natural extension of the classical
language containing the classical connectives (like modal languages are
extensions of the classical language). 

Furthermore, the new logic allows quantitative reasoning about amplitudes and probabilities, being in this respect much closer to the possible worlds logics for probability reasoning than to the traditional quantum logics. For other developments in this direction, also motivated by applications in quantum computation and information, see \cite{mey:pat:03a,mey:pat:03b}. 

Herein, we present a finitary Hilbert calculus for EQPL and show that it is weakly complete relatively to an oracle for arithmetical reasoning. Strong completeness is out of question since entailment is not compact. The proof of the weak completeness result was carried out using a non trivial extension of the technique proposed by Fagin, Halpern and Megiddo for simple probabilistic logics, together with three Henkin completions.

Although EQPL only provides the means for propositional, quantitative reasoning about quantum states, it is a mandatory step before further developments towards calculi for reasoning about the evolution of quantum systems (as already outlined in \cite{pmat:acs:04b}). The weak completeness result established here is interesting from the theoretical point of view and shows that the proposed language fits the proposed exogenous semantics. But, for practical applications in quantum system specification and verification, it seems better to go for model checking techniques.

Such future developments of our approach to quantum reasoning are briefly discussed in Section 6 of the paper. In Section 2, we briefly motivate the EQPL semantic concepts and key design ideas, directly based on the postulates of quantum physics. In Section 3, we present the EQPL language and semantics plus some examples. In Section 4, we introduce the axioms and rules of EQPL. Section 5 is fully dedicated to the proof of the main result (weak completeness of EQPL).

\section{Key design ideas}

Starting from the postulates of quantum mechanics (closely following \cite{coh:diu:lal:77}) we present the key ideas that guided the design of EQPL (together with a brief review of the relevant concepts and results of operator theory).

{\bf Postulate 1}: {\em Every isolated quantum system is described by a Hilbert space. The states of the quantum system are the unit vectors of the corresponding Hilbert space.}

Recall that a Hilbert space is a complete inner product space over $\comps$ (the field of complex numbers). For example, the states of an isolated qubit are vectors of the form
$z_0\ket{0}+z_1\ket{1}$ where $z_0,z_1\in\comps$ and $|z_0|^2+|z_1|^2=1$. In other words, they are unit vectors in the (unique up to isomorphism) Hilbert space of dimension two. Concerning EQPL, it is natural to represent each qubit by a propositional symbol (more appropriately called a qubit symbol). Furthermore, each qubit state (better called qubit valuation) should be a superposition of the two possible classical valuations.

{\bf Postulate 2}: {\em The Hilbert space of a quantum system composed of a
finite number of independent components is the tensor product of the
component Hilbert spaces.}

For example, 
$z_{00}\ket{00}+z_{01}\ket{01}+ z_{10}\ket{10}+z_{11}\ket{11}$, where $z_{00}, z_{10}, z_{01}, z_{11} \in \comps$ and $|z_{00}|^2+|z_{01}|^2+|z_{10}|^2+|z_{11}|^2=1$, is the general form of the states of an isolated pair of qubits. Returning to the design of EQPL, we conclude that we need two qubit symbols for working with two qubits. Moreover, in this case, a quantum valuation should be a superposition of the four possible classical valuations.

It is easy to generalize this idea to a finite set of qubits. However, as
usual in logic, we would like to work with a fixed, denumerable {\em alphabet of
qubit symbols}:
        $$\qB=\{\qb_k: k \in \nats\}.$$
But, then, what should be the Hilbert space for $\qB$? The answer, a key ingredient of the
envisaged EQPL semantics, is the Hilbert space $\Hil=\Hil(2^\qB)$ that we 
define by free construction from the set $2^\qB$ of all classical valuations over $\qB$. This free construction is as follows. Given an arbitrary set $V$, the Hilbert space $\Hil(V)$ is as follows:
\begin{itemize}
\item Each element of $\Hil$ is a map
$\ket{\psi}: V \to \comps$ such that:
\begin{itemize}
  \item $\supp(\ket{\psi})=\{v: \ket{\psi}(v)\neq0\}$ is countable;
  \item $\displaystyle
  \sum_{v\in V} |\ket{\psi}(v)|^2 =
  \sum_{v\in\supp(\ket{\psi})} |\ket{\psi}(v)|^2 < \infty$.
\end{itemize}

\item $\ket{\psi_1}+\ket{\psi_2} =
                \lambda v.\,\ket{\psi_1}(v)+\ket{\psi_2}(v)$.

\item $\alpha \ket{\psi} =
                \lambda v.\,\alpha \ket{\psi}(v)$.

\item $\bkt{\psi_1}{\psi_2} =
       \displaystyle\sum_{v\in V}\conjg{\ket{\psi_1}(v)}\,\ket{\psi_2}(v)$.

\end{itemize}
The inner product induces the norm
$\norm{\ket{\psi}}=\sqrt{\bkt{\psi}{\psi}}$ and, so, the distance
$d(\ket{\psi_1},\ket{\psi_2})=\norm{\ket{\psi_1}-\ket{\psi_2}}$. Since
$\Hil(V)$ is complete for this distance, $\Hil(V)$ is a Hilbert
space\footnote{Isomorphic to $L^2(V,\#)$ where $\#$ is the counting measure over $V$.}.

Given $v\in 2^\qB$, $\ket v$ is the vector of $\Hil$ defined as follows:
$\ket{v}(v)=1$ and $\ket{v}(v')=0$ for every $v'\neq v$. Observe that
$\{\ket{v}: v\in V\}$ is an orthonormal basis of $\Hil$. This basis
will play an important role in the semantics of EQPL and for this reason we
refer to it as being the {\em logical basis} of $\Hil$. 

The unit vectors of $\Hil$ are the envisaged quantum valuations over
$\qB$. Given a quantum valuation $\ket\psi$ and a classical valuation $v$, the inner product $\bkt{v}{\psi}$ is said to be the {\em logical amplitude} of $\ket\psi$ for $v$. As we shall see, these logical amplitudes are at the core of EQPL.

Observe that it is useful to be able to work with a constrained set $V$ of admissible classical valuations. That is, it is sometimes convenient to work with $V\subsetneq 2^\qb$. Indeed, we may want to impose classical constraints on the quantum valuations. For example, we may want to impose $(\qb_1 \ldisj \qb_2)$, constraining the quantum system to states giving amplitude zero to every valuation not satisfying this classical formula. Therefore, concerning the semantics of EQPL, we conclude that a quantum interpretation structure $\bfw$ should contain at least a set $V\subseteq 2^\qB$ (the set of admissible classical valuations) and a unit vector $\ket\psi$ in $\Hil$ (the quantum valuation or the quantum state) such that $\bkt{v}{\psi}=0$ for every $v \not\in V$.

Since we start with the semantics for the whole system (composed of the denumerable set $\qB$ of qubits), what is the role of Postulate 2? More precisely, how can we identify an independent subsystem? The solution is ``tensor factorization" that we proceed to explain.

Given $S \subseteq \qB$ and $V \subseteq 2^\qB$, we introduce 
$V_{[S]} = \{v|_{S}:v\in V \}$ and $V_{]S[} = \{v|_{\qB\setminus S}:v\in V \}$. We also need $\Hil_{[S]}=\Hil((2^\qB)_{[S]})$ and $\Hil_{]S[}=\Hil((2^\qB)_{]S[})$.
Then, $\Hil(V)$ is a subspace of $\Hil(2^\qB)$; 
$\Hil=\Hil_{[S]}\otimes\Hil_{]S[}$; and $\Hil(V) \subseteq \Hil(V_{[S]})\otimes\Hil(V_{]S[})$ where equality does not hold in general.

Given a unit $\ket{\psi}\in\Hil$, if there are unit $\ket{\psi'}\in\Hil_{[S]}$ and unit $\ket{\psi''}\in\Hil_{]S[}$ such that $\ket{\psi}=\ket{\psi'}\otimes\ket{\psi''}$ then we say that, at state $\ket{\psi}$, the qubits in $S$ are not entangled with those outside $S$. In this situation, the state $\ket\psi$ is said to be $S$-factorizable.
Furthermore, a vector $\ket{\psi} \in \Hil_{[S]}$ is said to be non factorizable if
there is no proper subset $S'$ of $S$ such that there are unit
$\ket{\psi'} \in \Hil_{[S']}$ and unit $\ket{\psi''} \in
\Hil_{[S\setminus S']}$ such that $\ket{\psi}=\ket{\psi'}\otimes\ket{\psi''}$.

Having in mind these semantic notions, given a finite set $F$ of qubit symbols, we conclude that the language of EQPL should provide the means for writing assertions about:
\begin{itemize}

\item non entanglement: ``the qubits in $F$ are not entangled with the other qubits" (that is, the quantum state at hand is $F$-factorizable); this assertion is made, as we shall see, with the EQPL formula $[F]$;

\item logical amplitudes: ``the amplitude of a classical valuation over $F$ is equal to a complex number"; that is, we need terms denoting arbitrary complex numbers and terms denoting logical amplitudes; more precisely, as we shall see, when the quantum state is $F$-factorizable, the EQPL term $\ket{\verum}_{FA}$ denotes the amplitude of the (unique) classical valuation $v^F_A$ over target $F$ that satisfies the qubits in $A \subseteq F$ and does not satisfy the qubits in $F \setminus A$.

\end{itemize}

Other useful quantum constructions will be introduced as abbreviations,
including {\it inter alia}:
\begin{itemize}

\item $[G|F]$ -- formula stating that the quantum state is $G$-factorizable if it is $F$-factorizable.

\item $\ket\alpha_{FA}$ -- term roughly denoting the amplitude of $v^F_A$ if this classical valuation satisfies $\alpha$, and equal to zero otherwise.

\item $\qiposs{F}{\alpha:u}$ -- formula stating that the quantum state is $F$-factorizable and that there is a classical valuation over $F$ in the $F$-component of the quantum state satisfying $\alpha$ and with non null amplitude $u$.

\end{itemize}

Unfortunately, the amplitude terms are not always meaningful on a given pair $\tuple{V,\ket\psi}$. Namely, they require that the target qubits are not entangled with the others. Therefore, we need more information in the envisaged notion of quantum interpretation structure. But, before we are ready to give the definition, we need some additional notation about partitions of $\qB$. Given a partition $\cS$ of $\qB$, let $\ucS$ be the set of all unions of elements of $\cS$. That is, $\ucS = \{\bigcup_{S\in\cR}S: \cR \subseteq \cS\}$.

A {\em quantum interpretation structure} is a tuple $\bfw=\tuple{V,\cS,\ket\psi,\nu}$ where:
\begin{itemize}

\item $V$ is a nonempty subset of $2^\qB$.

\item $\cS$ is a finite partition of $\qB$.

\item $\ket\psi=\{\ket\psi_{[R]}\}_{R\in\ucS}$ where each
$\ket\psi_{[R]}$ is a unit vector of $\Hil_{[R]}$ and such that:
\begin{enumerate}
  \item $\ket\psi_{[\emptyset]}=e^{i0}$;
  \item $\ket\psi_{[R]} =
        \displaystyle\bigotimes_{\scriptsize
        \begin{array}{l}S\in\cS\\S\subseteq R\end{array}}
        \ket\psi_{[S]}$
        for each non empty $R\in\ucS$;
  \item $\ket\psi_{[S]}$ is non factorizable for each $S\in\cS$;
  \item $\bkt{v}{\psi}_{[\qB]}=0$ if $v\not\in V$.
\end{enumerate}

\item $\nu:\{\nu_{FA}\}_{F\subseteq_\textrm{fin}\qB,A\subseteq F}$ where each $\nu_{FA}\in\comps$ and $\nu_{FA} = \bkt{v^F_A}{\psi}_{[F]}$ if $F \in \ucS$.

\end{itemize}

In such a structure we recognize the key elements $V$ (the set of admissible classical valuations) and $\ket\psi_{[\qB]}$ (the quantum state of the whole system). The additional information is the factorization of $\ket\psi_{[\qB]}$ and the map $\nu$ that provides the means for interpreting amplitude terms even when they are physically undefined. In this way we avoided the need to work with partial interpretation structures.
Observe also that, although we work in $\Hil=\Hil(2^\qB)$, clause 4 in the definition above imposes that (up to isomorphism) we only consider quantum states in $\Hil(V)$.

As we just saw, Postulates 1-2 were sufficient to guide us in the task of setting up the notion of quantum interpretation structure over which we shall be able to define the semantics of EQPL. Now, we turn our attention to the postulates concerning measurements of physical quantities.

{\bf Postulate 3}: {\em Every measurable physical quantity of an isolated quantum system is described by an observable acting on its Hilbert space.}

Recall that an observable is a Hermitian operator such that the direct sum of its eigensubspaces coincides with the underlying Hilbert space. Since the operator is Hermitian, its spectrum $\Omega$ (the set of its eigenvalues) is a subset of $\reals$. For each $e\in\Omega$, we denote the corresponding eigensubspace by $E_e$ and the projector onto $E_e$ by $P_e$.

{\bf Postulate 4}: {\em The possible outcomes of the measurement of a physical quantity are the eigenvalues of the corresponding observable. When the physical quantity is measured using observable $A$ on a system in a state $\ket{\psi}$, the resulting outcomes are ruled by the probability space %
$\Prob^A_\ket{\psi}=\tuple{\Omega,\cB|_\Omega,\mu^A_{\ket{\psi}}}$
where in the case of a countable spectrum
$$\mu^A_{\ket{\psi}}=\lambda B. \sum_{e\in\Omega} \chi_B(e)%
|P_e\ket{\psi}|^2\,.$$}

For the applications we have in mind in quantum computation and information, only {\it logical projective measurements over a finite set of qubits} are relevant. Given a quantum structure $\bfw=\tuple{V,\cS,\ket\psi,\nu}$, for each finite set $F$ of qubits, such measurements are defined using some observable $A_F$ on $\Hil$ such that:
\begin{itemize}

\item The spectrum of $A_F$ is equipotent\footnote{The chosen bijection depends on how the qubits are physically implemented. For example, when implementing a qubit using the spin of an electron, we may impose that spin $+\frac12$ corresponds to true and spin $-\frac12$ corresponds to false. But, as we shall see, the semantics of EQPL does not depend on the choice of the bijection, as long as one exists. The same happens in the case of classical logic -- its semantics does not depend on how bits are implemented. The details of which voltages correspond to which truth values are irrelevant.}  to $V_{[F]}$.

\item For each $v'\in V_{[F]}$, the corresponding eigenspace $E_{v'}$ is generated by
all vectors of the form $\ket{v'}\otimes \ket{v''}$ in $\Hil$. Thus,
each projector $P_{v'}$ is $\ket{v'}\bra{v'}\otimes 1_{\Hil_{]F[}}$.

\end{itemize}

For example, if the system is in the particular state 
$$\alpha_{00\omega_1}\ket{00\omega_1}+\alpha_{01\omega_2}\ket{01\omega_2}
+\alpha_{01\omega_3}\ket{01\omega_3}+
\alpha_{10\omega_4}\ket{10\omega_4}$$ 
then the probability of observing the first two qubits $\qb_0,\qb_1$ in the classical valuation $01$ is given by $|\alpha_{01\omega_2}|^2 + |\alpha_{01\omega_3}|^2$.

In general, the stochastic result of making a logical projective measurement of a finite set $F$ of qubits of the system at $\bfw=\tuple{V,\cS,\ket\psi,\nu}$ is fully described by the finite probability space $\Prob^F_\bfw = \tuple{V_{[F]},\wp{V_{[F]}},\mu^F_\bfw}$ where, for each $U\subseteq V_{[F]}$:
$$\mu^F_\bfw(U) = \sum_{v'\in U}\;\;%
\sum_{v''\in V_{]F[}}%
|\bkt{(v'\oplus v'')}{\psi}|^2\,.$$
Here, $v'\oplus v''$ denotes the (unique) classical valuation over all
qubits determined by $v'$ and $v''$.

Thus, we are able to say what is the probability in a given quantum state
of observing a classical formula $\alpha$ as being true. That is, given a quantum structure $\bfw$, we have the means for interpreting EQPL terms of the form $\prob\alpha$ that denote such probabilities.

Finally, although irrelevant to the design of EQPL, we mention {\it en passant} Postulate 5 that rules how quantum systems evolve.

{\bf Postulate 5}: {\em Excluding measurements, the evolution of a quantum
system is described by unitary transformations.}

This last postulate becomes relevant only when designing a dynamical extension of EQPL (see \cite{pmat:acs:04b}).

\section{Language and semantics}

The language of EQPL is composed of classical formulae, real terms, complex terms and quantum formulae that we proceed to introduce using an abstract version of the BNF notation \cite{nau:63} for a compact presentation of inductive definitions. For building terms, it is convenient to use real variables $X=\{x_k: k\in\nats\}$ and complex variables $Z=\{z_k: k\in\nats\}$.

Classical formulae: 
$$\alpha \;:=\; \qb \bnfd (\lneg\alpha) \bnfd (\alpha \limp \alpha)$$ 
As usual, other classical connectives like $\lconj,\ldisj,\leqv$, verum $\verum$ and falsum $\falsum$ are introduced as abbreviations. We denote the set of qubit symbols occurring in $\alpha$ by $\qB(\alpha)$. We say that a classical formula $\alpha$ is over a set $S$ of qubit symbols if $\qB(\alpha)\subseteq S$.

Real and complex terms (with the provisos computable real constant $r$, finite $F\subset\qB$ and $A \subseteq F$): 
$$\left\{\begin{array}{lll} t
&:=&\textstyle x \bnfd r \bnfd \prob\alpha \bnfd
    (t+t) \bnfd
    (t\, t) \bnfd \re(u) \bnfd \im(u) \bnfd
    \arg(u) \bnfd |u|\\
u &:=&z \bnfd \ket{\verum}_{FA} \bnfd (t+it) \bnfd
    t e^{it} \bnfd \conjg{u} \bnfd
    (u+u) \bnfd (u\,u) \bnfd \qif{\alpha}{u}{u}
\end{array}\right.$$
Most of these terms are self-explanatory or already motivated in the previous section. An explanation is needed concerning complex alternative terms: a term $\qif{\alpha}{u_1}{u_2}$ denotes the value denoted by $u_1$ if $\alpha$ is true, and denotes the value denoted by $u_2$ otherwise. 

Quantum formulae (with the proviso finite $F\subset\qB$):
$$\gamma \;:=\;  \alpha
  \bnfd (t \leq t) \bnfd [F] \bnfd
 (\qneg\gamma) \bnfd (\gamma\qimp\gamma)$$
Quantum negation $\qneg$ and quantum implication $\qimp$ are global
operators and should not be confused with their classical (local)
counterparts. As expected, other quantum connectives will be introduced
as abbreviations. But, before introducing the whole set of useful
abbreviations, we present the semantics of the language.

Given a set $S$ of qubit symbols and a set $V$ of valuations, the {\em
extent} at $V$ of classical formulae over $S$ is as follows:
\begin{itemize}
\item $\ext{\alpha}^S_V = \{v\in V_{[S]}: v\sat_c\alpha\}$.
\end{itemize}

By an {\em assignment} $\rho$ we mean a map such that $\rho(x)\in\reals$
for each $x\in X$ and $\rho(z)\in\comps$ for each $z\in Z$. 

The {\em denotation of terms} at $\bfw=\tuple{V,\cS,\ket\psi,\nu}$ and $\rho$ is inductively defined as follows (we refrain from spelling out the obvious clauses for interpreting arithmetical expressions): 
\begin{itemize}

\item $\termvalue{x}_{\bfw\rho} = \rho(x)$;

\item $\termvalue{r}_{\bfw\rho} = r$;

\item $\termvalue{\prob\alpha}_{\bfw\rho} =
       \mu^{\qB(\alpha)}_\bfw(\ext\alpha^{\qB(\alpha)}_V)$;

\item $\termvalue{z}_{\bfw\rho} = \rho(z)$;

\item $\termvalue{\ket\verum_{FA}}_{\bfw\rho} =
       \nu_{FA}$;

\item $\termvalue{\qif{\alpha}{u_1}{u_2}}_{\bfw\rho} =
       \left\{\begin{array}{ll}
       \termvalue{u_1}_{\bfw\rho} & \textrm{if } \bfw\rho\sat\alpha\\
       \termvalue{u_2}_{\bfw\rho} & \textrm{otherwise}
      \end{array}\right.$;
      
\item \dots

\end{itemize}

The {\em satisfaction of quantum formulae} at
$\bfw=\tuple{V,\cS,\ket\psi,\nu}$
and $\rho$ is inductively defined as follows:
\begin{itemize}

\item $\bfw\rho \sat \alpha$ iff
        $v \sat_c \alpha$ for every $v\in V$;

\item $\bfw\rho \sat (t_1 \leq t_2)$ iff
      $\termvalue{t_1}_{\bfw\rho} \leq \termvalue{t_2}_{\bfw\rho}$;

\item $\bfw\rho \sat [F]$ iff $F \in \ucS$;

\item $\bfw\rho \sat (\qneg\gamma)$ iff
        $\bfw\rho \not\sat \gamma$;

\item $\bfw\rho \sat (\gamma_1\qimp\gamma_2)$ iff
        $\bfw\rho \not\sat \gamma_1$ or $\bfw\rho \sat \gamma_2$.

\end{itemize}

As anticipated in the previous section, the proposed quantum language with the semantics above is rich enough to express interesting properties of quantum systems. To this end, it is quite useful to introduce other operations, connectives and modalities through
abbreviations. We start with some additional quantum connectives:

\begin{itemize}

\item quantum disjunction: $(\gamma_1 \qdisj \gamma_2)$ for
$((\qneg\gamma_1)\qimp\gamma_2)$;

\item quantum conjunction: $(\gamma_1 \qconj \gamma_2)$ for
$(\qneg((\qneg\gamma_1)\qdisj(\qneg\gamma_2)))$;

\item quantum equivalence: $(\gamma_1 \qeqv \gamma_2)$ for
$((\gamma_1\qimp\gamma_2)\qconj(\gamma_2\qimp\gamma_1))$.

\end{itemize}

Observe that the quantum connectives are classical in the sense that
quantum tautologies hold. For instance,
$(((\qneg\gamma_2)\qimp(\qneg\gamma_1))\qimp(\gamma_1\qimp\gamma_2))$ 
is satisfied by every quantum structure and assignment. But they do not coincide with the classical connectives! For instance, $(\lneg\alpha)$ entails $(\qneg\alpha)$ but not the other way around. For a more detailed discussion of the differences and relationship between these two versions of classical logic refer to \cite{pmat:acs:css:05}.

It is also useful to introduce some additional comparison predicate symbols:

\begin{itemize}

\item $(t_1<t_2)$ for $((t_1\leq t_2)\qconj
     (\qneg(t_2\leq t_1)))$;

\item $(t_1=t_2)$ for $((t_1\leq t_2)\qconj(t_2\leq t_1))$;

\item $(u_1=u_2)$ for $((\re(u_1)=\re(u_2))\qconj
(\im(u_1)=\im(u_2)))$.

\end{itemize}

Classical molecular formulae (classical conjunctions of literals) are used profusely in the sequel. To this end, we introduce the following abbreviation (with the provisos finite $F\subset\qB$ and $A\subseteq F$):
\begin{itemize}

\item $(\bigwedge_F A)$ for
$((\bigwedge_{\qb_k\in A}\;\qb_k) \lconj (\bigwedge_{\qb_k\in
(F\setminus A)}\,(\lneg\qb_k)))$.

\end{itemize}
Observe that the formula $(\bigwedge_F A)$ specifies the unique classical valuation $v^F_A$ over $F$ that satisfies the qubits in $A$ and does not satisfy the qubits in $F\setminus A$.

Logical amplitude terms are easily extended to any classical formula besides verum (with the provisos $\qB(\alpha)\subseteq F$, finite $F\subset\qB$ and $A\subseteq F$):

\begin{itemize}

\item $\ket\alpha_{FA}$ for
$\qif{((\wedge_F A)\limp\alpha)}{\ket\verum_{FA}}{0}$.

\end{itemize}
Intuitively, $\ket\alpha_{FA}$ coincides with $\ket\verum_{FA}$ if $v^F_A$ satisfies $\alpha$, and it is zero otherwise.

Logical amplitude vector terms are introduced as follows (with the proviso $\qB(\alpha)\subseteq F$):

\begin{itemize}

\item $\ket\alpha_F$ for $(\ket\alpha_{FA})_{A\subseteq F}$.

\end{itemize}

It turns out that it is convenient to introduce the additional syntactic category of logical amplitude vector terms for each finite set $F$ of qubit symbols:
$$\ket\omega_F = \ket\alpha_F \bnfd (u\,\ket\omega_F) \bnfd (\ket\omega_F
+ \ket\omega_F)$$
with the obvious abbreviation rules for multiplication by scalar and addition. Still concerning amplitude vector terms, the following abbreviations are handy:

\begin{itemize}

\item $\ket0_F$ for $(0\ket\verum_F)$;

\item $(\ket{\omega_1}_F=\ket{\omega_2}_F)$ for
$(\bigsqcap_{A\subseteq F} (\ket{\omega_1}_{FA}=\ket{\omega_2}_{FA}))$;

\item $(\ket{\omega_1}_F \subseteq \ket{\omega_2}_F)$ for
$(\bigsqcap_{A\subseteq F} ((\ket{\omega_1}_{FA}\neq0) \qimp
      (\ket{\omega_1}_{FA}=\ket{\omega_2}_{FA})))$.

\end{itemize}

Finally, we are ready to introduce the rest of the interesting quantum operations, predicates
and modalities:

\begin{itemize}

\item $[G|F]$ for
  $(\bigsqcap_{A'\subseteq G}\bigsqcap_{A''\subseteq F\setminus G}\;
  (\ket{\verum}_{F(A'A'')}= %
      \ket{\verum}_{GA'}\ket{\verum}_{(F\setminus G)A''}))$;

\item $(\qb_{k_1}\sim_F\qb_{k_2})$ for
      $(\qneg(\bigsqcup_{\scriptsize\begin{array}{l}
      G \subset F\\
      \qb_{k_1}\in G\\
      \qb_{k_2}\notin G\\
      \end{array}}[G]))$;

\item $\qiposs{F}{\alpha:u}$ for %
$([F] \qconj |u|>0 \qconj
   (\bigsqcup_{A\subseteq F} (\ket{\alpha}_{FA}=u)))$;

\item $\qiposs{F}{\alpha_1:u_1,\dots,\alpha_n:u_n}$ for
$(\qiposs{F}{\alpha_1:u_1} \qconj \dots \qconj
\qiposs{F}{\alpha_n:u_n})$;

\item $(\lposs \alpha)$ for $(0<\prob{\alpha})$;

\item $(\lnec \alpha)$ for $(1=\prob{\alpha})$.

\end{itemize}
Most of these quantum constructions were already discussed in
the previous section. The entanglement formula
$(\qb_{k_1}\sim_F\qb_{k_2})$ states that the two qubits are entangled.

Quantum molecular formulae (quantum conjunctions of literals) are also very useful. Note that a quantum literal is either a quantum atom or the quantum negation of a quantum atom. Looking at the grammar of quantum formulae, it is clear that quantum atoms are either classical formulae, or comparisons between real terms or non entanglement assertions: 
$$\qAtom := \alpha \bnfd (t \leq t) \bnfd [F]$$ 
To this end, we introduce the following abbreviation (with the provisos finite $Q \subset \qAtom$ and $D\subseteq Q$):
\begin{itemize}

\item $(\bigsqcap_Q D)$ for
$((\bigsqcap_{\delta\in D}\;\delta) \qconj (\bigsqcap_{\delta\in
(Q\setminus D)}\,(\qneg\delta)))$.

\end{itemize}
Observe that a quantum molecular formula defines a set of quantum structures that may be empty because, for instance, the quantum molecular formula $(\alpha\qconj(\lneg\alpha))$ has no models (here $Q=\{\alpha,(\lneg\alpha)\}=D$).

We finish this section with a simple example. Consider the following variant of Schrödinger's cat. The relevant attributes of the cat are: being inside or outside the box, alive or dead, and moving or still. These three attributes are represented by the qubits $\qb_0,\qb_1,\qb_2$, respectively. For the sake of readability we use instead $\cati,\cata,\catm$, respectively. The following EQPL formulae constrain the state of the cat at different levels of detail:

\begin{enumerate}

\item $[\cati,\cata,\catm]$;

\item $(\catm \limp \cata)$;

\item $((\lposs\,\cata) \qconj (\lposs\,(\lneg \cata)))$;

\item $(\qneg [\cata])$;

\item $(\prob\cata = \frac13)$;

\item $([\cata,\catm]\lposs\;
(\cata\lconj\catm):\frac{1}{\sqrt 6},\\[1mm] \hspace*{15mm}
(\cata\lconj(\lneg\catm)):\frac{1}{\sqrt 6},\\[1mm] \hspace*{15mm}
((\lneg\cata)\lconj(\lneg\catm)):e^{i\frac{\pi}{3}}\sqrt\frac{2}{3})$.

\end{enumerate}
Observe that the assertions above are consistent with each other. Intuitively, assertion 1 states that the qubits $\cati,\cata,\catm$ are not entangled with the other qubtis of the cat system. Assertion 2 is a classical constraint on the set of admissible valuations: if the cat is moving then it is alive. Assertion 3 states the famous paradox: the cat can be in a state where it is possible that the cat is alive and it is possible that the cat is dead. Assertion 4 states that the qubit $\cata$ is entangled with other qubits. Assertion 5 states that the cat is in a state where the probability of observing it alive (after collapsing the wave function) is $\frac13$. Finally, assertion 6 states that the qubits $\cata,\catm$ are not entangled with other qubits and that in the quantum state there is a classical valuation with amplitude $\frac1{\sqrt6}$ where the cat is alive and moving, there is another classical valuation also with amplitude $\frac1{\sqrt6}$ where the cat is alive and not moving, and there is a classical valuation with amplitude $e^{i\frac{\pi}{3}}\sqrt\frac{2}{3}$ where the cat is dead (and, thus, thanks to 2, also not moving). 

\section{Axiomatization}

Entailment for EQPL may be defined as expected -- we say that $\Gamma$ entails $\delta$, written $\Gamma\ent\eta$, if $\bfw\rho\sat\eta$ for every $\bfw$ and $\rho$ satisfying every element of $\Gamma$. 
But a finitely bounded version of entailment turns out to be more relevant. Given a finite set $F$ of qubit symbols, a quantum structure $\bfw=\tuple{V,\cS,\psi,\nu}$ is said to be $F$-factorizable if $F\in\ucS$.
Given a set $\Gamma$ of quantum formulae over $F$ and a quantum formula $\eta$ also over $F$, we say that the former $F$-entails the latter, written $\Gamma\ent_F\eta$ if $\bfw\rho\sat\eta$ for every $F$-factorizable $\bfw$ and $\rho$ satisfying every element of $\Gamma$. 

Observe that $\Gamma\ent\eta$ implies $\Gamma\ent_F\eta$ for every $F$. 
Furthermore, for any $\Gamma$ and $\eta$ over $F_1$, if $F_1\subseteq F_2$ and $\Gamma\ent_{F_2}\eta$ then $\Gamma\ent_{F_1}\eta$. 
Note also that $\Gamma,\eta_1\ent_F\eta_2$ iff $\Gamma\ent_F(\eta_1\qimp\eta_2)$, and a similar result holds for the unbounded entailment. That is, quantum implication does internalize the notion of quantum entailment in EQPL. 

It is also straightforward to verify that both entailments (unbounded and bounded) are not compact in the sense that there are $\Gamma$ and $\eta$ such that $\eta$ is entailed by $\Gamma$ but it is not entailed by any finite subset of $\Gamma$. Therefore, there is no hope of setting up a finitary axiomatization (that is, using only finitary rules) achieving strong completeness. But, it is possible to establish a finitary axiomatization that achieves $F$-{\em bounded weak completeness} for any finite $F$: $\ent_F \eta$ iff $\der_F \eta$. Indeed, the axioms and rules presented below are sound and adequate for $F$-validity as will be proved in the next section.

Before listing all axioms and rules we need to introduce the concept of {\em tautological quantum formula} or {\em quantum tautology}. A quantum formula $\gamma$ is said to be tautological if there are a classical tautology $\alpha$ and a substitution map $\sigma: \qB \to \qAtom$ such that $\gamma$ coincides with $\alpha^{\lneg,\limp}_{\qneg,\qimp}\sigma$. For instance, the quantum formula $((x_1 \leq x_2)\qimp(x_1 \leq x_2))$ is tautological (obtained, for example, from the classical tautology $(\qb_1\limp\qb_1)$). We also need the concept of arithmetical language:
$$\begin{array}{lll}
\upsilon &:=& (a \leq a) \bnfd
 (\qneg\upsilon) \bnfd (\upsilon\qimp\upsilon)\\
a &:=&\textstyle x \bnfd r \bnfd
    (a+a) \bnfd
    (a\, a) \bnfd \re(b) \bnfd \im(b) \bnfd
    \arg(b) \bnfd |b|\\
b &:=&z \bnfd (a+ia) \bnfd
    a e^{ia} \bnfd \conjg{b} \bnfd
    (b+b) \bnfd (b\, b)
\end{array}$$
Observe that an assignment $\rho$ is enough to interpret arithmetical formulae. An arithmetical formula $\upsilon$ is said to be valid if it is satisfied by every assignment. For instance, $(((t_1\leq t_2)\qconj(t_2\leq t_3)) \qimp(t_1\leq t_3))$ and $((u_1^2=-1)\qimp((u_1=i)\qdisj(u_1=-i)))$ are both universal arithmetical formulae (the latter using equality between complex numbers introduced as an abbreviation).

We are now ready to list the axioms and rules of our calculus for each finite set $F$ of qubit symbols:

\begin{itemize}
	
\item $\der_F \alpha$ for each classical tautology $\alpha$ \;\;[CTaut].

\item $\alpha_1,(\alpha_1\limp\alpha_2) \der_F \alpha_2$ \;\;[CMP].

\item $\der_F \gamma$ for each quantum tautology $\gamma$ \;\;[QTaut].

\item $\gamma_1,(\gamma_1\qimp\gamma_2) \der_F \gamma_2$ \;\;[QMP].

\item $\der_F \upsilon^{\vec x \vec z}_{\vec t \vec u}$ for each valid arithmetical formula $\upsilon$ \;\;[Oracle].

\item $\der_F ((\alpha_1\limp\alpha_2)\qimp (\alpha_1\qimp\alpha_2))$ \;\;[Lift$\limp$].

\item $\der_F ((\alpha_1\qconj\alpha_2)\qimp (\alpha_1\lconj\alpha_2))$ \;\;[Ref$\qconj$].

\item $\der_F (\alpha\qimp(\qif{\alpha}{u_1}{u_2}=u_1))$ \;\;[If$\verum$].

\item $\der_F ((\qneg\alpha)\qimp(\qif{\alpha}{u_1}{u_2}=u_2))$ \;\;[If$\falsum$].

\item $\der_F [F]$ \;\;[NEtg$F$].

\item $\der_F ([G_2] \qimp ([G_1] \qeqv [G_1|G_2]))$ for any $G_1\subseteq G_2$ \;\;[NEtg$|$].

\item $\der_F ([G_1] \qimp ([G_2] \qimp [G_1\cup G_2]))$ \;\;[NEtg$\cup$].

\item $\der_F ([G_1] \qimp ([G_2] \qimp [G_1\setminus G_2]))$ \;\;[NEtg$\setminus$].

\item $\der_F (\ket{\verum}_{\emptyset\emptyset}=1)$ \;\;[Empty].

\item $\der_F ((\lneg(\wedge_F A))\qimp(\ket\verum_{FA}=0))$ \;\;[NAdm].

\item $\der_F ([G]\qimp((\sum_{A\subseteq G}|\ket\verum_{GA}|^2)=1))$ \;\;[Unit].

\item $\der_F (\prob\alpha=(\sum_{A\subseteq F}|\ket\alpha_{FA}|^2))$ \;\;[Prob].

\end{itemize}

In total, we have only two rules ({\it modus ponens} for classical implication [CMP] and for quantum implication [QMP]\footnote{Actually, [CMP] can be derived from [QMP] and [Lift$\limp$].}) and fifteen axioms. The axioms are better understood in the following groups.

We have as axioms the classical tautologies and the quantum tautologies ([CTaut] and [QTaut], respectively). Since the set of classical tautologies and the set of quantum tautologies are both recursive, there is no need to spell out the details of tautological reasoning.

Axiom [Oracle] is more controversial -- we accept as an axiom any valid arithmetical formula. The set of valid arithmetical formulae is not even recursively enumerable, hence the name we chose for the axiom. We decided to use an arithmetical oracle for two reasons. First, we wanted to focus our attention on reasoning about quantum aspects without becoming lost in arithmetical details. And, second, the alternative of presenting a recursive axiomatization based on the theory of algebraic closed fields would require,  in order to maintain completeness, a relaxation of our semantics, maybe towards a point too far away from its intuitive roots in the postulates of quantum mechanics. However, this alternative is interesting also for other reasons and we shall come back to the issue in the last section of the paper.

Axioms [Lif$\limp$] and [Ref$\qconj$] are sufficient to relate (local) classical reasoning and (global) quantum tautological reasoning. Again, we refer to \cite{pmat:acs:css:05} for more details.

Axioms [If$\verum$] and [If$\falsum$] are self explanatory. They will be used in the completeness proof to remove alternative terms.

Axioms [NEtg$F$], [NEtg$|$], [NEtg$\cup$] and [NEtg$\setminus$] are enough to reason about non entanglement. Among other things they impose that non entanglement is closed under set theoretic operations (closure under intersection appears as a theorem as we shall see).

Axioms [Empty], [NAdm] and [Unit] rule logical amplitudes. Each of them closely reflects a property of our semantic structures.

Finally, axiom [Prob] relates probabilities and amplitudes, closely following Postulate 4.

As expected, we say that a formula $\eta$ over $F$ is an $F$-{\em theorem}, written $\der_F \eta$ if we can build a derivation of $\eta$ from the axioms using the rules (for $F$). As an illustration, consider the following derivation that establishes  for any finite $F$:
\begin{itemize}
	\item $\der_F (\prob\verum=1)$ \;\;[PUnit].
\end{itemize}
{\footnotesize\begin{tabbing}
\=0000\=0000000000000000000000000000000000000000000000000000000000000000\=\kill
\>1%
\>$[F]$%
\>NEtg$F$\\%
\ \\
\>2%
\>$([F]\qimp((\sum_{A\subseteq F}|\ket\verum_{FA}|^2)=1))$%
\>Unit\\%
\ \\
\>3%
\>$((\sum_{A\subseteq F}|\ket\verum_{FA}|^2)=1)$%
\>QMP:1,2\\%
\ \\
\>4%
\>$(\prob\verum=(\sum_{A\subseteq F}|\ket\verum_{FA}|^2))$%
\>Prob\\%
\ \\
\>5%
\>$((\prob\verum=(\sum_{A\subseteq F}|\ket\verum_{FA}|^2))\qimp
(((\sum_{A\subseteq F}|\ket\verum_{FA}|^2)=1)\qimp(\prob{\verum}=1)))$%
\>Oracle\\%
\ \\
\>6%
\>$(((\sum_{A\subseteq F}|\ket\verum_{FA}|^2)=1)\qimp(\prob{\verum}=1))$%
\>QMP:4,5\\%
\ \\
\>7%
\>$(\prob{\verum}=1)$%
\>QMP:3,6%
\end{tabbing}}

We finish this section with a list of interesting $F$-theorems. Some of them will be used in the proof of weak completeness presented in the next section, but others are mentioned just for illustration purposes.

The following theorem is a direct consequence of the non entanglement axioms and completes the picture of non entanglement being closed under set theoretic operations.
\begin{itemize}

\item $\der_F ([G_1] \qimp ([G_2] \qimp [G_1\cap G_2]))$ \;\;[NEtg$\cap$].

\end{itemize}
The following theorems give some insight on the major properties of logical amplitudes.
\begin{itemize}

\item $\der_F ((\ket{(\alpha_1\ldisj\alpha_2)}_G + \ket{(\alpha_1\lconj\alpha_2)}_G) = (\ket{\alpha_1}_G + \ket{\alpha_2}_G))$ \;\;[AAdd].

\item $\der_F ((\alpha_1\limp\alpha_2) \qimp
      (\ket{\alpha_1}_G \subseteq \ket{\alpha_2}_G))$ \;\;[AMon].

\item $\der_F ((\alpha_1 \leqv \alpha_2)\qimp
      (\ket{\alpha_1}_G = \ket{\alpha_2}_G))$ \;\;[ASoE].

\item $\der_F (\alpha \qimp (\ket{\alpha}_G = \ket{\verum}_G))$ \;\;[ANec].

\item $\der_F ((\ket\alpha_G +\ket{(\lneg\alpha)}_G) = \ket\verum_G)$ \;\;[AMExc].

\end{itemize}
The first of the following theorems about probability after measurements just states finite additivity. The second is an obvious instance of Postulate 4. The third relates logical reasoning with probability reasoning (monotonicity).
\begin{itemize}

\item $\der_F ((\prob{(\alpha_1\ldisj\alpha_2)} +
\prob{(\alpha_1 \lconj \alpha_2)}) = (\prob{\alpha_1} + \prob{\alpha_2}))$ \;\;[PAdd].

\item $\der_F (\qiposs{G}{(\bigwedge_G A):u}
      \qimp (\prob{(\bigwedge_G A)}=|u|^2))$ \;\;[Meas].

\item $\der_F ((\alpha_1\limp\alpha_2) \qimp
      (\prob{\alpha_1} \leq \prob{\alpha_2}))$ \;\;[PMon].

\end{itemize}
The following theorems show that the quantum and probability modalities do behave as normal modalities.
\begin{itemize}

\item $\der_F (\qiposs{G}{(\alpha\ldisj\alpha'):u}
             \qeqv
 (\qiposs{G}{\alpha:u} \qdisj \qiposs{G}{\alpha':u}))$ \;\;[QNorm].

\item $\der_F ((\alpha\limp\alpha') \qimp
(\qiposs{G}{\alpha:u} \qimp \qiposs{G}{\alpha':u}))$ \;\;[QMon].

\item $\der_F ((u=u') \qimp
(\qiposs{G}{\alpha:u} \qimp \qiposs{G}{\alpha:u'}))$ \;\;[QCong].

\item $\der_F (\alpha\qimp(\lnec\alpha))$ \;\;[PNec].

\item $\der_F ((\lnec(\alpha\limp\alpha'))\qimp((\lnec\alpha)\qimp(\lnec\alpha')))$ \;\;[PNorm].

\end{itemize}

\section{Proof of bounded weak completeness}

It is straightforward to prove that the calculus presented in the last section is strongly sound -- for any finite $F\subset\qB$, if $\Gamma \der_F \eta$ then $\Gamma \ent_F \eta$. Therefore, it is also weakly sound.

On the other hand, as already pointed out, it is not possible to achieve strong adequacy with a finitary calculus. But, for arbitrary finite $F\subset\qB$, we were able to prove $F$-bounded adequacy of the calculus -- if $\ent_F \eta$ then $\Gamma \der_F \eta$. Therefore, since we have soundness, our calculus is $F$-bounded weakly complete -- $\ent_F \eta$ iff $\der_F \eta$.

The quantitative nature of the language of EQPL raises specific problems when proving an adequacy result. These problems appear on top of those raised by the fact the calculus is not strongly complete. Thus, the traditional Henkin approach to adequacy proofs \cite{hen:50} is not the answer here, or, at least, is not the full answer. 

In the end, we were inspired by the Fagin-Halpern-Megiddo technique that was successfully applied in proving adequacy results for probability calculi \cite{fag:hal:meg:90}. The key step of this technique is the reduction of any formula to a disjunction of systems of linear inequations over the real numbers where each variable represents the probability of a classical molecular formula. A close exam of the technique suggests that it should be applicable (possibly after a suitable non trivial extension) to any quantitative logic where the disjunctive normal form lemma holds.

Actually, a quite significant revamp of the Fagin-Halpern-Megiddo technique was needed in order to cope with the novel aspects of EQPL: (i)~classical formulae mixed with arithmetic (in)equations; (ii)~global semantics of quantum connectives; (iii)~non entanglement atoms; (iv)~amplitude terms besides probability terms; and (v)~quantum structures instead of probability spaces. Note that the Fagin-Halpern-Meggido technique was first developed for a probabilistic logic somewhat simpler than the probabilistic fragment of EQPL.

In addition, we used the Henkin technique thrice: (i)~for removing alternative
terms; (ii)~for constructing the set of admissible valuations; and (iii)~for
building the finite partition of the set of qubits.

The rest of this section contains a step by step presentation of the proof of the $F$-bounded weak adequacy of EQPL.

Given a quantum formula $\gamma$ over $F$ we say that it is
$F$-consistent if $\not\der_F(\qneg\gamma)$. The proof is carried out by contraposition:
\begin{enumerate}
  \item Assume that $\not\der_F \gamma$.
  \item So, quantum tautologically, also $\not\der_F (\qneg(\qneg \gamma))$.
  \item Thus, $(\qneg\gamma)$ is $F$-consistent.
  \item Therefore, by the main lemma proved below,
        there are $F$-factorizable $\bfw$ and $\rho$
        such that $\bfw\rho\sat(\qneg\gamma)$.
  \item And, hence, it is not true that every such pair
  satisfies $\gamma$, that is, we established that $\not\ent_F\gamma$.
\end{enumerate}

It remains to prove the
model existence lemma: {\em If $\gamma$ is $F$-consistent then there are
$F$-factorizable $\bfw$ and $\rho$ such that $\bfw\rho\sat\gamma$.}

The quantum disjunctive normal form lemma holds in EQPL. Thus:
$$\der_F \left(\gamma \qeqv \bigsqcup_{D\in\qmols(\gamma)}
(\sqcap_{Q_\gamma} D)\right)$$\\%
where $\qmols(\gamma)= \{D\subseteq
Q_\gamma: \;\;\der_F ((\sqcap_{Q_\gamma} D) \qimp \gamma)\}$ and
$Q_\gamma$ is the set of $F$-quantum atoms used in $\gamma$.

Clearly, $\gamma$ is $F$-consistent iff
there is $D\in\qmols(\gamma)$ such that $(\sqcap_{Q_\gamma} D)$ is
$F$-consistent. Therefore, it is sufficient to prove the following
restricted model existence lemma: {\em If $(\sqcap_Q D)$ is
$F$-consistent then there are $F$-factor\-izable $\bfw$ and $\rho$ such
that $\bfw\rho \sat (\sqcap_Q D)$.}

Since $D = D_c \cup D_{\leq} \cup D_{[\,]}$, where
$D_c\subseteq Q_c = \{\alpha: \alpha \in Q\}$,
$D_{\leq}\subseteq Q_{\leq} = \{(t\leq t'): (t\leq t') \in Q\}$, and
$D_{[\,]}\subseteq Q_{[\,]} = \{[G]: [G] \in Q\}$, we have:
$$(\sqcap_Q D) =
((\sqcap_{Q_c} D_c) \sqcap (\sqcap_{Q_\leq} D_\leq) \sqcap
(\sqcap_{Q_{[\,]}} D_{[\,]})).$$

Our goal is to reduce everything to inequations. We start by getting rid of the non entanglement atoms.

Thanks to NEtg$F$ and NEtg$|$, we know that there is a quantum formula $\delta_{[]}$ without non entanglement atoms such that
$\der_F ((\sqcap_{Q_{[\,]}} D_{[\,]}) \qeqv \delta_{[\,]})$.
Thus,
$\der_F ((\sqcap_Q D) \qeqv \delta)$ where $\delta = ((\sqcap_{Q_c} D_c)
\sqcap (\sqcap_{Q_\leq} D_\leq) \sqcap \delta_{[\,]})$.

Note that $\delta_{[\,]}$ and, hence, $\delta$ are not
necessarily conjunctions of quantum literals (because it may happen that
a $[G]$ appears in $Q_{[\,]} \setminus D_{[\,]}$ and such a negation
involves a disjunction).
Using again the quantum disjunctive normal form lemma we have:
$$\der_F \left(\delta \qeqv \bigsqcup_{D\in\qmols(\delta)}
(\sqcap_{Q_\delta} D)\right).$$%
So, $\delta$ is $F$-consistent iff there is $D\in\qmols(\delta)$
such that $(\sqcap_{Q_\delta} D)$ is $F$-consistent.
Therefore, it is sufficient to prove the following
even more restricted model existence lemma: {\em If $(\sqcap_Q D)$
without entanglement atoms is $F$-consistent then there are
$F$-factorizable $\bfw$ and $\rho$ such that $\bfw\rho \sat (\sqcap_Q
D)$}.

Assume that $(\sqcap_Q D)$ is $F$-consistent and does not involve
non entanglement atoms (that is, $Q=Q_c\cup Q_\leq$ and $D=D_c\cup
D_\leq$). Our goal is to find an $F$-factorizable
$\bfw=\tuple{V,\cS,\ket\psi,\nu}$ and a $\rho$ satisfying this molecular
formula. We start by looking for $V$.

Before setting up $V$,
it is necessary to eliminate the probability and alternative terms and to
add maximally consistent information about the admissible classical
valuations. This desideratum is achieved as follows:

\begin{enumerate}

\item First, we replace  in $(\sqcap_Q D)$ each term $\prob\alpha$
by $\sum_{A\subseteq F}\ket\alpha_{FA}$.  Let %
$(\sqcap_{\overline{Q}} \overline D)$ be the result.

\item  Consider an ordering $\alpha_1,\dots,\alpha_m$ of the guards of
alternative terms occurring in $(\sqcap_{\overline{Q}} \overline D)$.

\item Consider the following sequence of formulae:
\begin{itemize}
  \item $\eta_0 = (\sqcap_{\overline{Q}} \overline D)$;
  \item $\eta_{k+1} =
  \left\{\begin{array}{ll}
  (\eta_k \qconj \alpha_k) & \textrm{if }
                     \der_F (\eta_k\qimp\alpha_k)\\
  (\eta_k \qconj (\qneg\alpha_k)) & \textrm{otherwise}
  \end{array}\right.$.
\end{itemize}

\item Observe that each $\eta_k$ is still $F$-consistent and a quantum
molecular formula. Furthermore, $\eta_m$ is maximal with respect to
guards.

\item Now we can replace each term $\qif{\alpha}{u_1}{u_2}$
occurring in $\eta_m$ by:
\begin{itemize}
\item $u_1$ if $\alpha$ is a quantum literal in $\eta_m$;
\item $u_2$ if $(\qneg\alpha)$ is a quantum literal in $\eta_m$.
\end{itemize}
Let $\overline\eta_m$ be the resulting formula.

\item  Consider an ordering $A_1,\dots,A_{m'}$ of the
subsets of $F$.

\item Consider the following sequence of formulae:
\begin{itemize}
  \item $\eta'_0 = \overline\eta_m$;
  \item $\eta'_{k+1} =
  \left\{\begin{array}{ll}
  (\eta'_k \qconj (\lneg(\wedge_F A_k))) & \textrm{if }
                     \der_F (\eta'_k\qimp(\lneg(\wedge_F A_k)))\\
  (\eta'_k \qconj (\qneg(\lneg(\wedge_F A_k)))) & \textrm{otherwise}
  \end{array}\right.$.
\end{itemize}

\item Observe that each $\eta'_k$ is still $F$-consistent and a quantum
molecular formula. Furthermore, $\eta'_{m'}$ does not contain probability
terms or alternative terms and is maximal with respect to admissible
classical valuations.

\item Thanks to Prob, If$\verum$ and If$\falsum$, denoting the resulting
still $F$-consistent molecular formula by $(\sqcap_{Q'} D') =
((\sqcap_{Q'_c} D'_c)\qconj(\sqcap_{Q'_\leq} D'_\leq))$, we have $$\der_F
((\sqcap_{Q'} D') \qimp (\sqcap_Q D)).$$

\item Therefore, we may proceed working towards the envisaged $\bfw$
and $\rho$ with the new formula.

\end{enumerate}

Having (while preserving $F$-consistency) eliminated
the probability and alternative terms and having determined the classical
valuations, we are ready to build $V$.
Let $V$ be composed of each $v\in2^\qB_{[F]}$
such that $v \sat \alpha$ for each $\alpha\in D'_c$.
Now we have to analyze two cases:

\begin{description}

\item[a)] Either for each $\alpha \in Q'_c \setminus D'_c$
there is a $v\in V$ such that $v \not\sat_c \alpha$ and, therefore,
this $V$ is viable because $$\tuple{V,\dots} \sat_F
(\sqcap_{Q'_c}D'_c)\,.$$

\item[b)] Or that is not the case. But, then, we would be able to contradict the $F$-consistency of $(\sqcap_{Q'}D')$ as follows:

\begin{enumerate}

\item Indeed, if it is not the case then there is a $\alpha\in
Q'_c\setminus D'_c$ such that $v \sat_c \alpha$ for all $v\in V$. That is, by construction of $V$, there is $\alpha\in Q'_c\setminus D'_c$
such that%
$$\ent_c \left(\left(\bigwedge_{\alpha'\in D'_c}\alpha'\right)\limp
\alpha\right).$$ %

\item So, by CTaut, there is $\alpha\in Q'_c\setminus D'_c$
such that%
$$\der_F \left(\left(\bigwedge_{\alpha'\in D'_c}\alpha'\right)\limp
\alpha\right).$$ %

\item Thus, by Lift$\limp$, there is $\alpha\in Q'_c\setminus D'_c$
such that%
$$\der_F \left(\left(\bigwedge_{\alpha'\in
D'_c}\alpha'\right)\qimp \alpha\right).$$ %

\item Thus, by Ref$\qconj$ and QTaut (transitivity of $\qimp$), there is
$\alpha\in Q'_c\setminus D'_c$
such that%
$$\der_F \left(\left(\bigsqcap_{\alpha'\in
D'_c}\alpha'\right)\qimp \alpha\right).$$ %

\item Therefore, by QTaut (right weakening of $\qimp$)%
$$\der_F \left(\left(\bigsqcap_{\alpha'\in D'_c}\alpha'\right)\qimp
\left(\bigsqcup_{\alpha\in Q'_c\setminus D'_c}\alpha\right)\right)$$ %
leading to %
$$\der_F \left(\qneg\left(\left(\bigsqcap_{\alpha'\in D'_c}\alpha'\right)\qconj
\left(\bigsqcap_{\alpha\in Q'_c\setminus D'_c}(\qneg\alpha)\right)\right)\right)$$ %
by several obvious tautological steps.

\item That is, we have $\der_F (\qneg(\sqcap_{Q'_c}D'_c))$,
contradicting the $F$-consistency of $(\sqcap_{Q'_c}D'_c)$.

\end{enumerate}

\end{description}

In short, we did find $V$ satisfying the classical part of $(\sqcap_{Q'} D')$. Let us proceed with the construction of the partition $\cS$.
The idea is to find a maximally fine partition $\cS_F$ of $F$
such that $((\sqcap_{Q'} D')\qconj (\sqcap_{S\in \cS_F} [S|F]))$ is
$F$-consistent, as follows:

\begin{enumerate}

\item Let $G_1,\dots,G_n$ be an ordering of the subsets of $F$.

\item Consider the following sequence of formulae:
\begin{itemize}
  \item $\gamma_0 = (\sqcap_Q' D')$;
  \item $\gamma_{k+1} =
  \left\{\begin{array}{ll}
  (\gamma_k \qconj [G_{k+1}|F]) & \textrm{if this formula is
                 $F$-consistent}\\
  \gamma_k & \textrm{otherwise}
  \end{array}\right.$.
\end{itemize}

\item Observe that each $\gamma_k$ is $F$-consistent  and, furthermore,
$\gamma_n$ is maximally so with respect to non entanglement assertions.

\item Let $U = \{G: [G|F] \textrm{ is a factor of } \gamma_n\}$.

\item Let $\cS_F$ be composed of all minimal (with respect to inclusion)
elements of $U$. Then, thanks to NEnt$\cap$, NEnt$\cup$ and
NEnt$\setminus$, it is straightforward to prove that $\cS_F$ is a
partition of $F$. Moreover, $\ucS_F = U$.

\item Let $\cS = \cS_F \cup \{\qB\setminus F\}$.
Observe that $\cS$ is finite.

\item Since $\der_F (\gamma_n \qimp \eta_m)$, we proceed working with
$\gamma_n$ in our task of completing the construction of $\bfw$ and
$\rho$.

\end{enumerate}

It remains to find $F$-factorizable $\ket\psi$, together with $\nu$ and $\rho$.
As already mentioned, the key idea is to reduce everything to a system of (in)equations on variables representing amplitudes. But, first we need to add the constraints imposed by the relevant axioms. Thanks to Unit, for every $G\in\ucS_F$, we can establish: 
$\der_F (\gamma_n \qimp ((\sum_{A \subseteq G}|\ket\verum_{GA}|^2) = 1))$.
Thanks to NAdm, for every $(\lneg(\wedge_F A))$ occurring in
$\gamma_n$, we have:
$\der_F (\gamma_n \qimp (\ket\verum_{FA} = 0))$.

Let ${\gamma_n^\bullet}$ be the formula
$$\left(\gamma_n \qconj %
\left(\bigsqcap_{G\in\ucS_F} %
\left(\left(\sum_{A \subseteq G}|\ket\verum_{GA}|^2\right) = 1\right)
\right) \qconj
\left(\bigsqcap_{(\lneg(\wedge_F A))\textrm{ in }\gamma_n} %
\left(\ket\verum_{FA} = 0\right) \right) \right).$$

Observe that we can derive: $\der_F (\gamma_n \qeqv {\gamma_n^\bullet})$.
Let $({\gamma_n^\bullet})_\leq$ the conjunction of the (in)equations in ${\gamma_n^\bullet}$. Consider the finite system of (in)equations obtained from
$({\gamma_n^\bullet})_\leq$ by replacing at each term of the form
$\ket\verum_{GA}$ by a fresh variable $z_{\ket\verum_{GA}}$.
Now we have to analyse two cases:

\begin{description}

\item[a)] Either the system of (in)equations has no solution.
But, in this case we would be able to contradict the $F$-consistency
of $\gamma_n$ as follows (using the arithmetical oracle):

\begin{enumerate}

\item Let $\Lambda_\leq$ be the (finite) set of arithmetic literals
occurring in $( {\gamma_n^\bullet})_\leq$ and $\Lambda_c$ be the (finite)
set of non-arithmetic literals in $( {\gamma_n^\bullet})_c$.

\item Since $( {\gamma_n^\bullet})_\leq = (\sqcap_{\upsilon \in
\Lambda_\leq}\upsilon)$, there is a bijection between $\Lambda_\leq$ and
the set of inequations composing the system described above.

\item From the fact that the system of inequations induced by $(
{\gamma_n^\bullet})_\leq$ has no solution, we conclude that there
is no assignment $\rho$ such that %
$\rho\sat \upsilon$ for all $\upsilon \in \Lambda_\leq$.

\item In other words, for all assignment $\rho$ there exists $\upsilon \in
\Lambda_\leq$ such that $\rho\sat(\qneg\upsilon)$ and so, thanks to
Oracle, we have: 
$\der_F (\sqcup_{\upsilon \in \Lambda_\leq}(\qneg\upsilon))$.

\item Hence, a fortiori, we obtain:   $\der_F ((\sqcup_{\gamma \in
\Lambda_c}(\qneg\gamma))\sqcup (\sqcup_{\upsilon \in
\Lambda_\leq}(\qneg\upsilon))).$

\item That is, since %
$$\begin{array}{rcl} ((\sqcup_{\gamma \in \Lambda_c}(\qneg\gamma))\sqcup
(\sqcup_{\upsilon \in
     \Lambda_\leq}(\qneg\upsilon))) &
     \qeqv & (\qneg((\sqcap_{\gamma \in \Lambda_c}\gamma)
          \sqcap (\sqcap_{\upsilon \in \Lambda_\leq}\upsilon)))\\[3mm]
& = & (\qneg(( {\gamma_n^\bullet})_c \sqcap
            ( {\gamma_n^\bullet})_\leq))\\[3mm]
& = & (\qneg  {\gamma_n^\bullet})\\[3mm]
& \qeqv & (\qneg \gamma_n),
\end{array}$$
we can conclude $\der_F (\qneg \gamma_n)$, contradicting the
$F$-consistency of $\gamma_n$.

\end{enumerate}

\item[b)] Or the system has at least one solution and
then we can build the envisaged $F$-factorizable $\bfw=\tuple{V,\cS,\ket\psi,\nu}$
and $\rho$ from any of the solutions in the following way:

\begin{itemize}

\item $V$ is as described above.

\item $\cS$ is as described above.

\item $\ket\psi = \{\ket\psi_{[R]}\}_{R\in\ucS}$ is obtained
as follows:

\begin{itemize}

\item $\ket\psi_{[G]}(v^G_A)$ is the solution value of
   $z_{\ket\verum_{GA}}$ for every $G \in
   \cS_F$ (note that $\ket\psi_{[G]}$
   is non-factorizable by construction of $\cS_F$);

\item $\ket\psi_{[\qB\setminus F]}$ is any non-factorizable unit vector in
      $\Hil(2^\qB_{[\qB\setminus F]})$ such that
      $\bkt{v}{\psi}_{[\qB\setminus F]}=0$ for every
      $v\not\in V_{[\qB\setminus F]}$;

\item $\ket\psi_{[\emptyset]}=e^{i0}$ and
      $\ket\psi_{[R]} =
        \otimes_{\scriptsize\begin{array}{l}S\in\cS\\S\subseteq R\end{array}}
        \ket\psi_{[S]}$
        for each non empty $R\in\ucS$.

\end{itemize}

\item $\nu=\{\nu_{GA}\}_{G \subset_\textrm{fin}\qB,A\subseteq G}$ is
chosen as follows:

\begin{itemize}

\item If $z_{\ket\verum_{GA}}$ is a variable of the system then $\nu_{GA}$
takes the value of this variable in the adopted solution.

\item Otherwise:

If $G \in \cS_F$ then %
$\nu_{GA} = \bkt{v^G_A}{\psi}_{[G]}$;

otherwise, the value of $\nu_{GA}$ can be chosen freely in
$\comps$.

\end{itemize}

\item $\rho$ is established as follows:
\begin{itemize}
\item $\rho(x)$ is equal to the value
      of $x$ if this variable occurs in the system,
      and given an arbitrary value otherwise;
\item $\rho(z)$ is equal to the value
      of $z$ if this variable occurs in the system,
      and given an arbitrary value otherwise.
\end{itemize}

\end{itemize}

Such a pair $\bfw\rho$ satisfies $({\gamma_n^\bullet})_\leq$ and, so,
also satisfies $(\sqcap_Q D)$.\hfill{QED}

\end{description}

\section{Concluding remarks}

Using a non trivial extension of the Fagin-Halpern-Megiddo technique together with three Henkin like completions we were able to prove the finitely bounded weak completeness of the proposed finitary axiomatization for EQPL. The arithmetical oracle was used once for obtaining a contradiction in the case where the induced system of (in)equations has no solution. 

The adoption of an arithmetical oracle for abstracting away the reasoning about real and complex numbers allowed us to concentrate on the quantum aspects of the calculus. Since the set of valid arithmetical formulae is not recursively enumerable there is no hope to find a recursive axiomatization while preserving weak completeness over the proposed semantics. But, it is viable and possible interesting to relax the semantics and replace the oracle with a recursive axiomatization of the algebraic closed fields. Parallel developments in probabilistic logic \cite{fag:hal:meg:90,aba:hal:94}, give us hope of obtaining even decidable calculi. But, then, we have to pay the price of working with relaxed quantum structures that are far away from their roots in the postulates of quantum mechanics. Nevertheless, this seems to be the solution towards model checking techniques for EQPL and its dynamical extensions. Such model checking techniques also require the development of the theory of quantum automata \cite{amm:pmat:acs:05}.

The weak completeness result obtained in this paper shows that the proposed language of EQPL is appropriate for the proposed exogenous semantics. Therefore, EQPL constitutes a sound basis for further developments of our approach to quantum reasoning, namely towards dynamical extensions for reasoning about the evolution of quantum systems and protocols. For preliminary results in this direction, see \cite{pmat:acs:04b} where DEQPL (a dynamical extension of EQPL) is outlined. Recent work on dynamical versions of traditional quantum logic \cite{bal:sme:04} should also be taken into account.
Another interesting development, also from the applications point of view, will be directed at a EQFOL (a FOL version of exogenous quantum logic).

The detailed analysis of the weak completeness proof reinforces the idea (already present in the choice of the EQPL abbreviations) of the key role, when using EQPL for reasoning, of a finite context of qubit symbols. One wonders if this assumption can be relaxed to any recursive set of qubits by starting with classical $\omega$-infinitary propositional logic \cite{kar:64}. At least from a theoretical point of view, this line of work should be explored.

As we saw, the semantics of EQPL is based on pure quantum states of collections of qubits. Recall that pure quantum states are unit vectors of the underlying Hilbert space. In consequence, EQPL provides the means for asserting properties of and reason about such vectors. Therefore, EQPL is not insensitive to the global phase of the quantum state. One may argue that it should be insensitive since no physical measurement will ever be able to distinguish two quantum states that are equivalent up to global phase. We decided to make EQPL as it is (that is, sensitive to global phase) for two reasons. In practice, physicists and quantum computer scientists need to work with both levels of abstraction. Sometimes they want to work with states as unit vectors. Sometimes they want to abstract away the global phase. Therefore, a calculus supporting the former level of abstraction is also useful. The second reason is a consequence of the fact that forgetting global phase requires a major semantic shift. Indeed, it is better solved by identifying a quantum state, not with a unit vector of the underlying Hilbert space, but, instead, with a density operator working on that space, that is, working with ensembles or mixed quantum states in general. 

Such shift towards a semantics based on density operators will lead to a quite different quantum logic (but still extending classical logic by applying the exogenous approach) that will also be useful for reasoning about quantum systems evolving under partial tracing, besides unitary transformations and measurements. Clearly, this is yet another line of research that will deserve attention.

Finally, the relationship between the exogenous quantum logics and the more traditional quantum logics (based on the original Birkhoff and von Neumann proposal) should be explored. At the preliminary stage of work in this direction, it seems that most of the qualitative assertions possible in the latter can be made in the former and that most of the quantitative assertions possible in the former can be borrowed by extensions of the latter.


\section*{Acknowledgments}

The authors wish to express their deep gratitude to the regular participants in the QCI Seminar at CLC, specially Ana Maria Martins and Vítor Rocha Vieira, who attended early presentations of EQPL and gave very useful feedback that helped us to get over our initial difficulties in and misunderstandings of quantum physics. This work was partially supported by FCT and FEDER through POCTI, namely
via QuantLog POCTI/MAT/55796/2004 Project.

\end{document}